\documentclass[10pt,a4paper,twoside]{article}
\usepackage{amsmath,amssymb,enumerate,
,theorem}
\usepackage{epsfig,color,graphics,graphicx}
\usepackage{amssymb,amsbsy,amsmath,amsfonts,amssymb,amscd}
\usepackage[english]{babel}
\newcommand{\R}{{\mathbb{R}}}

\newcommand{\assign}{:=}

\newcommand{\tmem}[1]{{\em #1\/}}
\newcommand{\tmop}[1]{\ensuremath{\operatorname{#1}}}
\newcommand{\tmtextbf}[1]{{\bfseries{#1}}}
\newcommand{\tmtextit}[1]{{\itshape{#1}}}
\newcommand{\tmtextsc}[1]{{\scshape{#1}}}

\newenvironment{itemizedot}{\begin{itemize}  }{\end{itemize}}
\newenvironment{proof}{\noindent\textbf{Proof\ }}{\hspace*{\fill}$\Box$\medskip}
\newtheorem{lemma}{Lemma}[section]
\newtheorem{theorem}{Theorem}
\newtheorem{proposition}[lemma]{Proposition}

\theoremstyle{definition}

\newtheorem{definitionproposition}[lemma]{Definition - Proposition}
\newtheorem{remark}{Remark}

\begin{document}

\title{On the cost of null-control of an artificial advection-diffusion problem.}
\author{Pierre Cornilleau\thanks{Teacher at Lyc\'ee du parc des Loges, 1, boulevard des Champs-\'Elys\'ees, 91012 \'Evry, France.\newline e-mail: pierre.cornilleau@ens-lyon.org.}  \& Sergio Guerrero \thanks{Laboratoire Jacques-Louis Lions,
Universit\'e Pierre et Marie Curie, 75252 Paris C\'edex 05, France. \newline e-mail : guerrero@ann.jussieu.fr.}}\maketitle

\begin{abstract}
In this paper we study the null-controllability of an artificial advection-diffusion system in dimension $n$.
Using a spectral method, we prove that the control cost goes to zero exponentially when the viscosity vanishes and the control time is large enough. On the other hand,  we prove that the control cost tends to infinity exponentially when the viscosity vanishes and the control time is small enough.
\end{abstract}

\tableofcontents

\section*{Introduction}
 The following paper continues \cite{cG} and deals with an advection-diffusion problem with small viscosity
truncated in one space direction. This problem was first considered in \cite{H}, where the Cauchy problem has been studied when the viscosity tends to zero.

\paragraph{Artificial advection-diffusion problem}
In this paper, we consider an advection-diffusion system in a strip $\Omega:=\{(x',x_n)\in\mathbb{R}^{n-1} \times (-L,0)\}$ ($n\ge1$ and $L$ a positive constant) with particular artificial boundary conditions on both sides of the domain. As indicated above, this system was considered in \cite{H} (see section 6 in that reference):
\begin{equation}\label{systemeHalpern} \left\{ \begin{array}{c}
     u_t + \partial_{x_n} u - \varepsilon \Delta u = 0\\
     \varepsilon(u_t + \partial_{\nu} u) = 0\\
     \varepsilon (u_t + \partial_{\nu} u) + u = 0\\
     u (0, .) = u_0
   \end{array} \right. \begin{array}{c}
     \text{in } (0, T) \times \Omega,\\
     \text{on } (0, T) \times \Gamma_0,\\
     \text{on } (0, T) \times \Gamma_1,\\
     \text{in } \Omega,
   \end{array} \end{equation}
where $T>0$, $\Gamma_0 \assign \mathbb{R}^{n - 1} \times \{0\}$, $\Gamma_1
\assign \mathbb{R}^{n - 1} \times \{- L\}$ and we have denoted $\partial_{x_n}$ the partial derivative with respect to $x_n$ and $\partial_{\nu}$ the normal derivative.

We are here interested in the uniform boundary controllability of (\ref{systemeHalpern}):
    $$
    (S_v)\left\{ \begin{array}{c}
     u_t + \partial_{x_n} u - \varepsilon \Delta u = 0\\

     \varepsilon (u_t + \partial_{\nu} u) + u1_{\Gamma_1} = v1_{\Gamma_0}\\
     u (0, .) = u_0
   \end{array} \right. \begin{array}{c}
     \text{in } (0, T) \times \Omega,\\
     \text{on } (0, T) \times \partial\Omega,\\
     \text{in } \Omega.
   \end{array}
   $$
   We recall that, if $X$ is defined as the closure of $\mathcal{C}^{\infty} (
\bar{\Omega}$) for the norm
$$
\|u\|_X \assign \left( \|u\|^2_{L^2(\Omega)} +
\varepsilon \|u\|^2_{L^2(\partial\Omega)} \right)^{\frac{1}{2}},
$$
the system $(S_v)$ is well-posed in this space (see section 1 in \cite{cG}).

  In this present paper, we study  the so-called \textit{null controllability} of this system on $\Gamma_0$
\[
\tmop{for} \tmop{given} u_0\in X, \tmop{find} v\in L^2((0,T),\Gamma_0) \tmop{such} \tmop{that}
   \tmop{the} \tmop{solution} \tmop{of} (S_v) \tmop{satisfies} u (T) \equiv 0.
   \]
Furthermore, we will be interested in the continuous dependence of these controls on the initial data, that is to say, the existence of $C>0$ such that
\begin{equation}\label{OI}
 \Vert v \Vert_{L^2((0,T),\Gamma_0)} \leq C \Vert u_0\Vert_{X}, \quad \forall u_0 \in X.
\end{equation}
We will denote by $C (\varepsilon)$ the cost of the null-control, which is the
smallest constant $C$ fulfilling estimate \eqref{OI}. We remark that $C(\varepsilon)$ equals $+\infty$ when the null controllability does not hold.

\par \vskip0.2cm
We have proved in \cite{cG} that $(S_v)$ is null-controllable in dimension $n=1$ for any $T,\,\varepsilon>0$. However, the argument of Miller \cite{M} cannot be directly applied in this situation (for more details see the appendix in \cite{cG}).

\par \vskip0.2cm

In the present paper, we first obtain a precise upper bound on the null-control cost using a spectral approach combined with a Carleman estimate in dimension one.
In a second part, we use a more classical method to prove that for $T$ small enough, the cost $C(\varepsilon)$ exponentially tends to infinity when $\varepsilon \to 0$.

\par\vskip0.2cm

In the context of degeneration of a parabolic-to-hyperbolic type systems, similar results have been obtained by many authors in dimension one (see, for instance, \cite{CG,G} (one dimensional heat equation) and \cite{GG} (linear Korteweg de Vries equation))
but also in dimension $n$ (see \cite{GL}).
However, our results seem to be new in the context  of a system which lacks of regularizing effect. A reasonable conjecture seems to be that the system is not null controllable for small $T, \varepsilon>0$.

\paragraph{Main results}

Our main results are the following:

\begin{theorem}\label{controln}
If $T / L$ is large enough, the cost of the null-control $C (\varepsilon)$ tends to zero exponentially as
  $\varepsilon \rightarrow 0$:
$$\exists C,k>0 \text{ such that } \quad C(\varepsilon)\le C e^{-k/\varepsilon} \quad \forall \varepsilon\in(0,1).$$
\end{theorem}

\begin{remark} \vskip0.2cm $ $
\begin{itemize}
\item One can in fact obtain the same controllability result when the control acts on $\Gamma_1$ (see also Remark 4).
\item
  The fact that the control cost tends to zero tells intuitively that the
  state almost vanishes for $T / L$ big enough. This is to be connected
  with the fact that, for $\varepsilon = 0$, the system is purely advective
  and then that, for $T > L$, its state vanishes.
  \end{itemize}
\end{remark}

\begin{theorem}\label{noncontroln}
If $T < L$, the cost of the null-control $C(\varepsilon)$ exponentially tends to infinity when $\varepsilon \to 0$:
$$
\forall T<L, \ \exists \varepsilon_0>0, \ \exists C,k>0 \text{ such that }\ \forall \varepsilon \in (0, \varepsilon_0), \ C(\varepsilon) \ge C e^{k/\varepsilon}.
$$
\end{theorem}

\begin{remark} This result is analogous to other results already obtained in the context of vanishing viscosity
(see for instance \cite[Theorem 2]{CG}, \cite[Theorem 1]{GL}). Observe that in these papers, the null controllability for small $\varepsilon$ and $T$ was known while in the present situation this question is open. 
\end{remark}

As usual in the context of linear controllability problems, we introduce the following adjoint system:
\[ (S') \left\{ \begin{array}{c}
     \varphi_t + \partial_{x_n} \varphi + \varepsilon \Delta \varphi = 0\\
     \varepsilon (\varphi_t - \partial_{\nu} \varphi) - \varphi = 0\\
     \varphi_t - \partial_{\nu} \varphi = 0\\
     \varphi (T, .) = \varphi_T
   \end{array} \right. \begin{array}{c}
     \text{in } (0, T) \times \Omega,\\
     \text{on } (0, T) \times \Gamma_0,\\
     \text{on } (0, T) \times \Gamma_1,\\
     \text{in } \Omega,
   \end{array} \]
where $\varphi_T \in X$. It is classical to prove that the controllability of system $(S_v)$ and the observability of system $(S')$
are equivalent (see, for instance, \cite{DR}):
\begin{proposition}\label{ObsControln}
  The following properties are equivalent :
  \begin{itemizedot}
    \item $\exists C_1>0,  \forall \varphi_T \in X_{}$; $\|\varphi (0,.)
    \|_X \leq C_1 \| \varphi \|_{L^2 ((0, T),\Gamma_0)}$ where
    $\varphi$ is the solution of problem $(S'),$

    \item $\exists C_2>0,  \forall u_0 \in X, \exists v \in L^2((0, T), \Gamma_0)$ such that $\| v \|_{L^2( (0, T), \Gamma_0) } \leqslant C_2 \| u_0 \|_X$ and \
    the solution $u$ of $(S_v)$ satisfies $u (T) = 0$.
  \end{itemizedot}
  Moreover, $C_1 = C_2$.
\end{proposition}

The rest of the article is organized as follows: in the first section, we introduce a one-dimensional problem with parameter and study its well-posedness and its observability. For the latter, we show a Carleman inequality for the associated adjoint system (see Proposition 1.4 below).
We consequently deduce Theorem \ref{controln} and Theorem \ref{noncontroln} in section two. In the appendix, we furthermore give a proof of Proposition 1.4.

Moreover, we note that the substitution $(t,x_n) \rightarrow (L t, Lx_n)$ allows us to assume that $L=1$. This hypothesis will be imposed until the end of the paper.

\par $\ $

\par \textbf{Notations:}
\par \noindent $A \lesssim B$ means that, for some universal constant $c>0$, $A\leq c B.$
\par \noindent $A \sim B$ means that, for some universal constant $c>1$, $c^{-1} B \leq A \leq c B.$

\section{A one-dimensional problem with parameter}
In all this section, we assume that $n=1$.
We also denote $X^1$ the space $X$.
We shall prove the following null-controllability result:
\begin{proposition}\label{coutcontrole}
  If $T$ is sufficiently large, there exists $\varepsilon_0>0$ such that, for any  $a \ge 0$ and $\varepsilon \in (0, \varepsilon_0)$, the system
$$
    (S^a_v)\left\{ \begin{array}{c}
     u_t + u_x - \varepsilon u_{xx}+ au = 0\\

     \varepsilon (u_t + \partial_{\nu} u)= v\\
\varepsilon(u_t+ \partial_{\nu} u)+u=0\\
     u (0, .) = u_0
   \end{array} \right. \begin{array}{c}
     \text{in } (0, T) \times (-1,0),\\
     \text{on } (0, T) \times \{0\},\\
\text{on } (0, T) \times \{-1\},\\
     \text{in } (-1,0),
   \end{array}
   $$
is null-controllable. That is to say, for any $u_0 \in X^1$ there exists $v \in L^2(0,T)$ such that the solution $u$ of $(S^a_v)$ satisfies
$$
u(T)\equiv 0 \,\,\text{ and }\,\,\|v\|_{L^2(0,T)} \le C(\varepsilon,a) \|u_0\|_{X^1}.
$$
Moreover, the cost $C(\varepsilon,a)$ is bounded by
$$ C \exp \left( - \frac{k}{\varepsilon} \right)$$
  where $C$ and $k$ are some positive constants independent from $a$ and $\varepsilon$.
\end{proposition}

\subsection{Cauchy problem and duality}

First, we briefly show that problem $(S^a_v)$ is well-posed.
\\Indeed, we consider the
 bilinear form defined by
\begin{equation*}\label{b-def}
\forall u_1,u_2 \in H^1(-1,0), \quad b (u_1, u_2) = \varepsilon\int_{-1}^0 \partial_x u_1. \partial_x u_2+\int_{-1}^0  u_2 \partial_{x} u_1 + u_1(-1) u_2(-1) .
\end{equation*}
With the help of this bilinear form, one may now consider the space
\begin{equation*}
{\cal D}:=\left\{u_1\in X^1; \sup_{u_2 \in {\cal C}^{\infty}([-1,0]); \ ||u_2||_{X^1}\le 1} |b(u_1,u_2)|<+\infty \right\}
\end{equation*}
 equipped with the natural norm
$$\| u_1\| _{\mathcal{D}}=\| u_1\|_{X^1}+\sup_{u_2 \in {\cal C}^{\infty}([-1,0]); \ ||u_2||_{X^1}\le 1} |b(u_1,u_2)|.$$
Note that, using an integration by parts, one shows that $b(u_1,u_2)$ is well-defined for $u_1\in X^1$ and $u_2\in {\cal C}^\infty([-1,0])$ and that the map
$$ u_2 \in X\mapsto b(u_1,u_2) \in \mathbb{R}$$
is well-defined and continuous for any $u_1 \in {\cal D}$. Using the Riesz representation theorem, we can define a maximal monotone operator ${\cal A}$ with domain ${\cal D(A)}={\cal D}$ and such that
$$\forall u_1 \in {\cal D(A)}, \ \forall u_2\in X^1, \quad <-{\cal A}u_1, u_2>_{X^1}=b(u_1,u_2).$$
The Riesz representation theorem also provides the existence of a dissipative bounded operator ${\cal B}$ on $X^1$ such that
$$\forall u_1,u_2 \in X^1, \ <{\cal B}(u_1),u_2>_{X^1}=- \int_{-1}^0 u_1 u_2.$$
Using Rellich theorem, one easily sees that ${\cal B}$ is ${\cal A}$-compact (according to Definition 2.15 of \cite[Chapter III]{EN}) i.e. that
$$ {\cal B}: {\cal D(A)} \to X^1 \text{ is compact}$$
and, using Corollary 2.17 of \cite[Chapter III]{EN}, we get that the operator ${\cal A}+ a{\cal B}$ generates a contraction semi-group on $X^1$ for any $a \ge 0$.
Since $(S^a_0)$ can be written in the following
abstract way
\[ \left\{ \begin{array}{c}
    u_t = (\mathcal{A}+a {\cal B}) u,\\
     u(0, .) = u_0,
   \end{array} \right.
 \]
we have shown that the homogeneous problem $(S^a_0)$ possesses, for any $u_0 \in X^{1}$, a unique solution  $u \in \mathcal{C} ( [0,T], X^1)$. We will
call these solutions  \tmem{weak solutions} opposed to \tmem{strong
solutions} \ i.e. such that $u_0 \in \mathcal{D(A)}$ and which fulfill $
u \in \mathcal{C} ( \mathbb{R}^+, \mathcal{D(A)}) \cap \mathcal{C}^1 (\mathbb{R}^+, X^1)$.

\vspace{0.3 cm}

We now conclude as in Proposition 5 of \cite{cG} to the existence and uniqueness of solution to the nonhomogeneous problem $(S^a_v)$. More precisely, one has the following:
\begin{definitionproposition}\label{bienpose}$ $ \vspace{0.1cm}
\begin{itemize}
\item For $f\in L^2((0,T)\times(-1,0))$, $g_0\in L^2((0,T))$ and $g_1\in L^2((0,T))$, we put
 \[(S^a_{f,g_0,g_1}) \left\{ \begin{array}{c}
     u_t + u_x - \varepsilon u_{xx}+ a u = f\\
     \varepsilon(u_t + \partial_{\nu} u) = g_0\\
     \varepsilon (u_t + \partial_{\nu} u) + u = g_1\\
     u (0, .) = u_0
   \end{array} \right. \begin{array}{c}
     \text{in } (0, T) \times (-1,0),\\
     \text{on } (0, T) \times \{ 0\},\\
     \text{on } (0, T) \times \{ -1\},\\
     \text{in } (-1,0),
   \end{array}  \]
 and we say that $u \in \mathcal{C} ( \left[ 0, T \right], X^1)$ is a solution of
  $(S^a_{f,g_0,g_1})$ if, for every function $\psi \in \mathcal{C} ([0, T], {\cal D(A^*)}) \cap \mathcal{C}^1 ([0, T ], X^1)$, the following
  identity holds :
  \[ \int_0^{\tau} \left( <u, \psi_t>_{X^1}+<u,({\cal A}+a {\cal B})^* \psi>_{X^1} + <F,\psi>_{X^1}\right) = \left[ <u(t),\psi(t)>_{X^1}\right]_{t=0}^{t=\tau}
  \quad \forall \tau\in [0,T],
  \]
where we have defined, using the Riesz representation theorem, $F(t)\in X^1$ such that
\begin{equation*}\label{defF}
<F(t),\phi>_{X^1}=\int_{-1}^0 f(t)\phi+\int_{\{-1,0\}} g(t)\phi, \quad \forall \phi\in X^1.
\end{equation*}
and $g$ is a function on $(0,T) \times \{-1,0\}$ such that $g=g_0$ on $(0,T) \times \{0\}$, $g=g_1$ on $(0,T) \times \{-1\}$.
\item Let $T > 0$, $u_0 \in X$, $f\in L^2((0,T)\times(-1,0))$, $g_0\in L^2((0,T))$ and $g_1\in L^2((0,T))$. Then
  $(S^a_{f,g_0,g_1})$ possesses a unique solution $u$.
\end{itemize}
\end{definitionproposition}
\begin{proof}
This proof being very similar to the one of Proposition 5 of \cite{cG}, we think that a sketch will suffice.
\\First, if $u$ belongs to  $u\in \mathcal{C} ([0, T], {\cal D(A)}) \cap \mathcal{C}^1 ([0, T ], X^1)$ then,
using Duhamel formula and the density of ${\cal D}({\cal A}^*)$ in $X^1$, one obtains that
$u$ is a solution of $(S^a_{f,g_0,g_1})$ if and only if
\[  u(t) = e^{t ({\cal A}+a{\cal B}) }u_0 + \int_0^t e^{(t - s) ({\cal A}+a{\cal B}) } F (s)ds, \quad \forall t\in[0,T].\]
The general case now follow by a standard approximation argument.
\end{proof}

In order to study the null-controllability of system $(S^a_v)$, we shall focus on its adjoint problem, namely:
\[ (S'^a) \left\{ \begin{array}{c}
     \varphi_t + \varphi_x + \varepsilon \varphi_{xx} -a \varphi = 0\\
     \varepsilon (\varphi_t - \partial_{\nu} \varphi) - \varphi = 0\\
     \varphi_t - \partial_{\nu} \varphi = 0\\
     \varphi (T, .) = \varphi_T
   \end{array} \right. \begin{array}{c}
     \text{in } (0, T) \times (-1,0),\\
     \text{on } (0, T) \times \{0\},\\
     \text{on } (0, T) \times \{-1\},\\
     \text{in } (-1,0).
   \end{array} \]
An analogous semigroup method as presented above show that the adjoint problem $(S'^a)$ possesses, for any $\varphi_T \in X^{1}$, a unique solution  $\varphi \in \mathcal{C} ( [0,T], X^1)$ such that
\begin{equation}\label{dissipative}\forall t \in [0,T], \ \|\varphi(t)\|_{X^1} \le \|\varphi_T\|_{X^1}.\end{equation}
\begin{remark}\label{IDK}
This estimate also holds for solutions to system $(S')$. Indeed, the associated operator generates a contraction semigroup on $X$ (see \cite[Section 1.1]{cG}).
\end{remark}

In the following proposition, we also recall without proof the classical equivalence between observability and controllability.
\begin{proposition}\label{ObsControl}
  The following properties are equivalent :
  \begin{itemizedot}
    \item $\exists C_1>0,  \forall \varphi_T \in X^{1}$; $\|\varphi (0)
    \|_{X^1} \leq C_1 \| \varphi(.,0) \|_{L^2 (0, T)}$ where
    $\varphi$ is the solution of problem $(S'^a),$

    \item $\exists C_2>0,  \forall u_0 \in X^1, \exists v \in L^2(0, T)$ such that $\| v \|_{L^2
    (0, T) } \leqslant C_2 \| u_0 \|_{X^1}$ and \
    the solution $u$ of problem $(S^a_v)$satisfies $u (T) = 0$.
  \end{itemizedot}
  Moreover, $C_1 = C_2$.
\end{proposition}

\subsection{Proof of Proposition \ref{coutcontrole}}
\subsubsection{Carleman inequality}

In this paragraph, we state a Carleman-type inequality keeping track of the explicit dependence of all the constants with respect to $a$, $\varepsilon$ and $T$. As in \cite{FI}, we introduce the following weight functions:
\[ \forall x\in [-1,0], \quad \eta(x) \assign 2 + x,\,\,\,\, \alpha(t,x) \assign \frac{e^3 - e^{\eta(x)}}{t ( T - t)}, \,\,\,\,\phi(t,x)
   \assign \frac{e^{ \eta(x)}}{t ({T} - t)} . \]
\par \noindent One may show the following Carleman inequality.
\begin{proposition}\label{The:Carleman}
There exists $C>0$ and $s_0>0$ such that for every $\varepsilon\in (0,1)$ and every $s \geqslant s_0 ( \varepsilon^{-1} (T + T^2)+a^{1/2} \varepsilon^{-1/2} T^2)$ the following inequality is satisfied for every $\varphi_T\in X$:
\begin{equation}\label{Carlemana}
\begin{array}{l}\displaystyle
 s^3 \int_{(0,T)\times (-1,0)} \phi^3
     e^{- 2 s \alpha} | \varphi |^2
  + s^3 \int_{(0, T)\times \{0,-1\}} \phi^3 e^{- 2 s \alpha} | \varphi
     |^2
\leqslant C s^7\int_{(0, T)\times \{0\}}e^{-4s\alpha+2s\alpha(.,-1)}\phi^7|\varphi|^2.
     \end{array}
\end{equation}
Here, $\varphi$ stands for the solution of $(S'^a)$ associated to $\varphi_T$.
\end{proposition}

This Carleman estimate is quite similar to the one obtained in \cite[Theorem 9]{cG}. We have thus postponed its proof to appendix A.
\begin{remark}\label{Betis}
One can in fact obtain the following Carleman estimate with control term in $\Gamma_1$ :
\begin{equation}\label{Espagne}
\begin{array}{l}\displaystyle
 s^3 \int_{(0,T)\times (-1,0)}\phi^3
     e^{- 2 s \alpha} | \varphi |^2 + s^3 \int_{(0, T)\times \{0,-1\}} \phi^3 e^{- 2 s \alpha} | \varphi
     |^2
\leqslant C s^7\int_{(0, T)\times \{-1\}}e^{-4s\alpha+2s\alpha(.,0)}\phi^7|\varphi|^2,
     \end{array}
\end{equation}
simply by choosing the weight function $\eta(x)$ equal to $x\mapsto -x+1$ - the proof being very similar. This inequality is the first ingredient to prove
the first point stated in Remark 1.
\end{remark}

\subsubsection{Dissipation result}

In this paragraph, we show a dissipation result for the solutions of $(S'^a)$. We will distinguish two cases depending on the size of $a$.
\begin{itemize}

 \item \underline{Case $a\le \varepsilon^{-1}$}.

  Inspired by \cite{Dan}, we introduce a weight function $\theta (x) = \exp(\frac{\lambda}{\varepsilon} x)$ for some constant $\lambda \in(0,1)$ which will be fixed below.

  We multiply the first equation in $(S'^a)$ by $\theta \varphi$
  and we integrate on $(-1,0)$. This gives :
  \[
  \frac{1}{2} \frac{d}{dt} \left( \int_{-1}^0 \theta | \varphi
     |^2 \right) =  \underbrace{-\int_{-1}^0 \theta \varphi \varphi_x- \varepsilon \int_{-1}^0
     \theta \varphi \varphi_{xx}}_{A}+a\int_{-1}^0\theta|\varphi|^2 .
     \]
  Using now $\theta' =  \frac{\lambda}{\varepsilon} \theta $ and integrating by parts several times, we obtain
    \begin{eqnarray*}
  A&=&\frac{\lambda}{2\varepsilon}(1-\lambda)\int_{-1}^0\theta |\varphi|^2+\varepsilon\int_{-1}^0\theta |\varphi_x|^2 + \frac{1-\lambda}{2}\left( -\theta(0)
     |\varphi(\cdot,0)|^2 + \theta(-1) |\varphi(\cdot,-1)|^2
     \right)
     \\
     \\
     &-& \varepsilon \left(\theta(0)
     \varphi(\cdot,0) \varphi_x(\cdot,0)-\theta(-1) \varphi(\cdot,-1)\varphi_x(\cdot,-1)
     \right) .
     \end{eqnarray*}
  Using now the boundary conditions for $\varphi$ (see $(S'^a)$) and the fact that $a\geq 0$,
  we get
  \[
  \frac{d}{dt} \left( \int_{-1}^0 \theta | \varphi |^2 \right)+2 \varepsilon \int_{\{-1,0\}} \theta \varphi_t \varphi \geq
     \frac{\lambda(1- \lambda)}{\varepsilon} \int_{-1}^0 \theta | \varphi |^2 + (1-\lambda)
     \int_{\{-1,0\}} \theta | \varphi |^2  .
     \]
  Since $\lambda \in (0,1)$, we readily deduce
  \[ \frac{d}{dt} \left( \| \sqrt{\theta} \varphi (t)\|^2_{X^1} \right) \geq
     \frac{\lambda(1- \lambda)}{\varepsilon} \| \sqrt{\theta} \varphi
     (t)\|^2_{X^1}. \]
Gronwall's lemma combined with $\exp(-\frac{\lambda}{\varepsilon}) \leqslant \theta \leqslant 1$ successively gives, for $0 \leq t_1 \leq t_2 \leq T$,
  \[
  \| \sqrt{\theta} \varphi(t_1)\|^2_{X^1} \leqslant \exp \left( -\frac{\lambda(1- \lambda)}{\varepsilon} (t_2 - t_1)\right) \| \sqrt{\theta} \varphi (t_2)\|^2_{X^1} \]
  and
  \[ \| \varphi (t_1)\|^2_{X^1} \leqslant \exp \left( - \frac{1}{\varepsilon}\left( \lambda(1- \lambda)(t_2 - t_1)-\lambda \right)  \right) \| \varphi (t_2)\|^2_{X^1} \]
  For $t_2-t_1>1$, we finally choose
  $$
  \lambda:=\frac{t_2-t_1-1}{2(t_2-t_1)} \in(0,1),
  $$
  which gives
   \[
   \| \varphi (t_1)\|_{X^1} \leqslant \exp\left\{- \frac{(t_2 - t_1 -
     1)^2}{4\varepsilon (t_2-t_1)}\right\} \| \varphi (t_2)\|_{X^1}
     \]
     if $t_2-t_1>1$.

\item \underline{Case $a\ge \varepsilon^{-1} $}.

We multiply the equation satisfied by $\varphi$ by $ \varphi$  and we integrate on $(-1,0)$. We get the following identity, after an integration by parts in space:
  \begin{eqnarray*} \frac{1}{2} \frac{d}{dt} \left( \int_{-1}^0  | \varphi
     |^2 \right) &=&  -\frac{1}{2} \int_{-1}^0  \partial_{x}(| \varphi |^2)- \varepsilon \int_{-1}^0\varphi_{xx} \varphi+ a\int_{-1}^0 | \varphi
     |^2  .\\
&=& -\frac{1}{2}( |\varphi(\cdot,0)|^2-|\varphi(\cdot,-1)|^2)-\varepsilon  \varphi_x(\cdot,0)\varphi(\cdot,0) +\varepsilon  \varphi_x(\cdot,-1)\varphi(\cdot,-1)\\
&+&\varepsilon \int_{-1}^0 |\varphi_x|^2 +a\int_{-1}^0 | \varphi|^2.
\end{eqnarray*}
 Using now the boundary conditions,  we easily deduce
  \[
  \frac{d}{dt} \| \varphi(\cdot) \|_{X^1}^2  =|\varphi(\cdot,0)|^2+|\varphi(\cdot,-1)|^2+2\varepsilon \int_{-1}^0 |\varphi_x|^2 +2a\int_{-1}^0 | \varphi|^2.
  \]
  On the other hand, a standard trace result gives, for some constant $c\in ]0,1]$ (see for instance \cite[Theorem 1.5.10]{Gr})
  \[
  c a^{1/2} \varepsilon^{1/2}(|\varphi(\cdot,0)|^2+|\varphi(\cdot,-1)|^2) \le \varepsilon \int_{-1}^0 |\varphi_x|^2 +a\int_{-1}^0 | \varphi|^2 \]
  and, consequently,  we get, using that $a \ge \varepsilon^{-1}$,
  \[
  \frac{d}{dt} \left( \| \varphi (\cdot)\|^2_{X^1} \right) \geq
     ca^{1/2} \varepsilon^{-1/2} \| \varphi
     (\cdot)\|^2_{X^1}.
     \]
  Gronwall's lemma finally gives, for any $0 \le t_1 \le t_2 \le T$,
  \begin{equation*}
 \|  \varphi(t_1)\|^2_{X^1} \leqslant \exp \left( -ca^{1/2} \varepsilon^{-1/2} (t_2 - t_1)\right) \| \varphi (t_2)\|^2_{X^1} .
  \end{equation*}

  Summing up, we have shown the following dissipation result:
  \begin{lemma}
   There exists $c_0>0$ such that, for any $\varepsilon\in (0,1)$, $a\geq 0$, $t_2-t_1>1$ and any solution $\varphi$ of $(S'^a)$,
    \begin{equation} \label{dissipation}
  \|  \varphi(t_1)\|^2_{X^1} \leqslant \exp \left( -c_0\max\{a^{1/2},\varepsilon^{-1/2}\} \varepsilon^{-1/2} \frac{(t_2 - t_1 -1)^2}{t_2-t_1}\right) \| \varphi (t_2)\|^2_{X^1}.
  \end{equation}
  \end{lemma}

\subsubsection{Observability result}

We estimate both sides of the Carleman inequality obtained in Proposition \ref{The:Carleman}. Putting $m = e^3 - e^{2}$ and $M
    = e^3-e$, we first have
    \[
    s^7\int_{(0, T)\times \{0\}}e^{-4s\alpha+2s\alpha(.,-1)}\phi^7|\varphi|^2 \lesssim s^7 T^{- 14} \exp \left( \frac{s(8M-16m)
       }{T^2} \right) \int_{(0,T) \times \{0\}} | \varphi |^2.
       \]
    On the other hand, using that $\displaystyle{\phi \gtrsim \frac{1}{ T^2}}$ on $\left. [
    \frac{T}{4}, \frac{3 T}{4} \right]$, we have the following estimate from below for the left hand-side of the Carleman inequality \eqref{Carlemana}
    \[
    \frac{s^3}{T^6} \exp
       \left( - \frac{32sM}{ T^2} \right) \left(\int_{\frac{T}{4}}^{\frac{3
       T}{4}} \int_{- 1}^0 | \varphi |^2 + \int_{\frac{T}{4}}^{\frac{3 T}{4}}
       \int_{\left. \{- 1, 0 \right\}} | \varphi |^2\right).
       \]
   Consequently we get that
    \[
    \| \varphi \|^2_{L^2((T/4,3T/4); X^1)}
       \lesssim C \int_{(0,T) \times \{0\}} | \varphi |^2,
       \]
where $C=s^4 T^{- 8} \exp \left( \frac{16s(3M-m)
       }{ T^2} \right)$.   Choosing now $s
\sim \varepsilon^{-1}( T + \max\{(a\varepsilon)^{1/2},1\}  T^2)$, $C$ is estimated by, for some $c'>0$ independent from $T\ge1$,
    \[
    \varepsilon^{-4}\max\{(a\varepsilon)^2,1\} \exp\left(c'  \varepsilon^{-1}\max\{(a\varepsilon)^{1/2},1\}\right) \lesssim \exp\left(c'' \varepsilon^{-1}\max\{(a\varepsilon)^{1/2},1\}\right)
    \]
    for any $c'' > c'$. Summing up, we have obtained
    \begin{equation}\label{estimationcout1}
     \| \varphi \|^2_{L^2((T/4,3T/4); X^1)} \lesssim \exp\left(c'' \varepsilon^{-1}\max\{(a\varepsilon)^{1/2},1\}\right) \int_{(0,T) \times \{0\}}
       | \varphi |^2 .
    \end{equation}

    We now use the dissipation property \eqref{dissipation} with $t_1 = 0$ and $ t_2 = t \in
    \left] \frac{T}{4}, \frac{3 T}{4} \right[$. We easily get, for $T\ge 8$,
    \begin{equation}\label{estimationcout2} \frac{T}{2} \exp\left(\frac{c_0 T}{16} \varepsilon^{-1}\max\{(a\varepsilon)^{1/2},1\} \right) \| \varphi (0) \|^2_{X^1} \leqslant \| \varphi \|^2_{L^2((T/4,3T/4); X^1)}.
    \end{equation}
Combining \eqref{estimationcout1} with \eqref{estimationcout2} finally gives the result with moreover
$$k = \frac{c_0T}{16}-c''>0 \Longleftrightarrow T>16\frac{c''}{c_0},$$
using Proposition \ref{ObsControl}.
\end{itemize}

\section{Proof of the main results}

We are now able to deduce Theorem \ref{controln} and Theorem \ref{noncontroln}.
\vspace{0.3cm}
\\   As long as Theorem \ref{controln} is concerned, we will show that the cost associated to the null controllability problem
$$
    (S_v)\left\{ \begin{array}{c}
     u_t + \partial_{x_n} u - \varepsilon \Delta u = 0\\
\varepsilon (u_t + \partial_{\nu} u)= v\\
\varepsilon(u_t+ \partial_{\nu} u)+u=0\\
     u (0, .) = u_0
   \end{array} \right. \begin{array}{c}
     \text{in } (0, T) \times \Omega,\\
     \text{on } (0, T) \times \Gamma_0,\\
  \text{on } (0, T) \times \Gamma_1,\\
     \text{in } \Omega,
   \end{array}
   $$
can be estimated using a Fourier transform in $x'$.
\\ On the other hand, we will use a standard approach combining a dissipation result and a kind of conservation of energy to prove Theorem \ref{noncontroln}.
\\We first define, for any $f \in X$ and for a.e. $\xi' \in \R^{n-1}$, the Fourier transform of $f$ with respect to $x'$ by
$$ \hat{f}^{\xi'}(x_n)=\int_{\R^{n-1}} e^{- i\xi'.x'} f(x',x_n) d x'.$$
For real-valued functions $f$, we also define its real and imaginary part by, for a.e. $\xi' \in \R^{n-1}$,
$$ \hat{f}_r^{\xi'}(x_n)=\int_{\R^{n-1}} \cos( \xi'.x') f(x',x_n) d x' \ \text{ and } \ \hat{f}_i^{\xi'}( x_n)=-\int_{\R^{n-1}} \sin( \xi'.x') f(x',x_n) d x'.$$

\subsection{Proof of Theorem \ref{controln}}
We make use of Proposition \ref{coutcontrole}. We obtain that, for $T$ sufficiently large, $\varepsilon$ sufficiently small and for a.e. $\xi'\in \R^{n-1}$, there exists $v_r^{\xi'} \in L^2(0,T)$ such that
the solution $\hat{u}_r^{\xi'}$ of
$$
   \left\{ \begin{array}{c}
     \partial_t\hat{u}_r^{\xi'} + \partial_{x_n} \hat{u}_r^{\xi'} - \varepsilon \partial^2_{x_n}\hat{u}_r^{\xi'}+ \varepsilon |\xi'|^2 \hat{u}_r^{\xi'} = 0\\

     \varepsilon (\partial_t \hat{u}_r^{\xi'} + \partial_{\nu} \hat{u}_r^{\xi'})= v_r^{\xi'}\\
\varepsilon(\partial_t \hat{u}_r^{\xi'}+ \partial_{\nu} \hat{u}_r^{\xi'})+\hat{u}_r^{\xi'}=0\\
     \hat{u}_r^{\xi'} (0, .) = \hat{u_0}_r^{\xi'}
   \end{array} \right. \begin{array}{c}
     \text{in } (0, T) \times (-1,0),\\
     \text{on } (0, T) \times \{0\},\\
\text{on } (0, T) \times \{-1\},\\
     \text{in } (-1,0),
   \end{array}
   $$
satisfies
$$  \hat{u}_r^{\xi'}(T)\equiv 0$$
and
$$ \left\| v_r^{\xi'}\right\|_{L^2(0,T)} \le C \exp\left( -\frac{k}{\varepsilon}\right) \left\|\hat{u_0}_r^{\xi'}\right\|_{X^1}.$$
Using analogous notations for the imaginary part, we deduce that, putting $v^{\xi'}= v_r^{\xi'}-i v_i^{\xi'}$, the solution
 of
$$
   \left\{ \begin{array}{c}
     \partial_t\hat{u}^{\xi'} + \partial_{x_n} \hat{u}^{\xi'} - \varepsilon \partial^2_{x_n}\hat{u}^{\xi'}+ \varepsilon |\xi'|^2 \hat{u}^{\xi'} = 0\\

     \varepsilon (\partial_t \hat{u}^{\xi'} + \partial_{\nu} \hat{u}^{\xi'})= v^{\xi'}\\
\varepsilon(\partial_t \hat{u}^{\xi'}+ \partial_{\nu} \hat{u}^{\xi'})+\hat{u}^{\xi'}=0\\
     \hat{u}^{\xi'} (0, .) = \hat{u_0}^{\xi'}
   \end{array} \right. \begin{array}{c}
     \text{in } (0, T) \times (-1,0),\\
     \text{on } (0, T) \times \{0\},\\
\text{on } (0, T) \times \{-1\},\\
     \text{in } (-1,0),
   \end{array}
   $$
satisfies
$$  \hat{u}^{\xi'}(T)\equiv 0$$
and
$$ \left\| v^{\xi'}\right\|_{L^2(0,T)} \le C \exp\left(- \frac{k}{\varepsilon}\right) \left\|\hat{u_0}^{\xi'}\right\|_{X^1}.$$
It is now straightforward that, defining $v$ as the inverse Fourier transform of $\xi' \mapsto v^{\xi'}$, the solution of
$$
    (S_v)\left\{ \begin{array}{c}
     u_t + \partial_{x_n} u - \varepsilon \Delta u = 0\\
\varepsilon (u_t + \partial_{\nu} u)= v\\
\varepsilon(u_t+ \partial_{\nu} u)+u=0\\
     u (0, .) = u_0
   \end{array} \right. \begin{array}{c}
     \text{in } (0, T) \times \Omega,\\
     \text{on } (0, T) \times \Gamma_0,\\
  \text{on } (0, T) \times \Gamma_1,\\
     \text{in } \Omega,
   \end{array}
   $$
satisfies
$$  u(T)\equiv 0$$
and, using Parseval-Plancherel's identity,
$$ \left\| v\right\|_{L^2((0,T),\Gamma_0)} \le C \exp\left( -\frac{k}{\varepsilon}\right) \left\|u_0\right\|_{X}.$$
This ends the proof.

\subsection{Proof of Theorem \ref{noncontroln} }

In this paragraph, we follow the method exposed in \cite[Section 5]{GG} to get a lower bound on the cost of null-control. More precisely, we are going to find a function $\varphi_T$ such that the associated solution to $(S')$ satisfies
\begin{equation}\label{numberone}
\|\varphi\|_{L^2((0,T),\Gamma_0)}\lesssim e^{-C/\varepsilon}
\end{equation}
and
\begin{equation}\label{numbertwo}
\|\varphi(0,\cdot)\|_X\gtrsim 1,
\end{equation}
whenever $\varepsilon$ is small enough and $T<1$.

Let $\delta>0$ small enough such that $ 4\delta <1-T$  and let $\varphi_T$ be a smooth function defined in $\Omega=\R^{n-1}\times (-1,0)$ such that
\begin{equation}\label{supports} \left\{ \begin{array}{c} \text{Supp}(\varphi_T) \subset  \R^{n-1}\times (-2\delta, -\delta),\\
\displaystyle\|\varphi_T\|^2_X=\int_{\Omega} |\varphi_T|^2=1.\\
      \end{array} \right.
      \end{equation}

\begin{itemize}
 \item  Proof of (\ref{numberone}).

 We consider $\rho(x_n)=\exp\{\lambda\varepsilon^{-1}x_n\}$ for all $x_n\in (-1,0)$ and some $\lambda\in (0,1)$. Furthermore, we define a function $\Psi\in {\cal C}^{\infty}(\mathbb{R})$ such that
 $$
 \left\{\begin{array}{ll}
 \Psi=0&\hbox{ in }(-\infty,-3\delta),
 \\ \noalign{\medskip}
 \Psi=1&\hbox{ in }(-2\delta,+\infty),
  \\ \noalign{\medskip}
  \Psi'\geq 0
 \end{array}
 \right.
 $$
 and denote
 $$
 \psi_j(t,x):=\Psi^{(j)}(x_n+T-t)\quad 0\leq j\leq 2.
 $$

 Then, we multiply the equation in $(S')$ by $2\rho \psi_0 \varphi$ and we integrate in $\Omega$ :
 \begin{equation}\label{numberthree}
 \begin{array}{l}\displaystyle
 -\frac{1}{2}\frac{d}{dt}\int_{\Omega}\rho\psi_0|\varphi|^2=-\varepsilon\int_{\Omega}\rho\psi_0\Delta\varphi\varphi
 -\int_{\Omega}\rho\psi_0\partial_{x_n}\varphi\varphi-\frac{1}{2}\int_{\Omega}\rho\psi_1|\varphi|^2.
 \end{array}
 \end{equation}
 Integrating by parts in the first term of the right-hand side, we have
 $$
 \begin{array}{l}\displaystyle
 -\varepsilon\int_{\Omega}\rho\psi_0\Delta\varphi\varphi=-\varepsilon\int_{\Gamma}\rho\psi_0\partial_{\nu}\varphi\varphi
 +\lambda\int_{\Omega}\rho\psi_0\partial_{x_n}\varphi\varphi+\varepsilon\int_{\Omega}\rho\psi_1\partial_{x_n}\varphi\varphi+\varepsilon\int_{\Omega}\rho\psi_0|\nabla\varphi|^2.
\end{array}
$$
We use the boundary conditions in $(S')$ for the first term and we integrate by parts again in the second and third term. This yields :
$$
\begin{array}{l}\displaystyle
-\varepsilon\int_{\Omega}\rho\psi_0\Delta\varphi\varphi=-\frac{\varepsilon}{2}\frac{d}{dt}\int_{\Gamma}\rho\psi_0|\varphi|^2+\left(\frac{\lambda}{2}+1
\right)\int_{\Gamma_0}\rho\psi_0|\varphi|^2
\\ \noalign{\medskip}\displaystyle
\phantom{-\varepsilon\int_{\Omega}\rho\psi_0\Delta\varphi\varphi}-\varepsilon\int_{\Gamma_1}\rho\psi_1|\varphi|^2-\frac{\lambda}{2}\int_{\Gamma_1}\rho\psi_0|\varphi|^2+\varepsilon\int_{\Omega}\rho\psi_0|\nabla\varphi|^2
\\ \noalign{\medskip}\displaystyle
\phantom{-\varepsilon\int_{\Omega}\rho\psi_0\Delta\varphi\varphi}
-\frac{\lambda^2}{2\varepsilon}\int_{\Omega}\rho\psi_0|\varphi|^2-\lambda\int_{\Omega}\rho\psi_1|\varphi|^2-\frac{\varepsilon}{2}\int_{\Omega}\rho\psi_2|\varphi|^2.
\end{array}
$$

We plug this into (\ref{numberthree}) and we integrate by parts in the second term of the right-hand side of (\ref{numberthree}). We obtain :
$$
\begin{array}{l}\displaystyle
\frac{d}{dt}\int_{\Omega}\rho\psi_0|\varphi|^2=2\varepsilon\int_{\Omega}\rho\psi_0|\nabla\varphi|^2+\frac{\lambda(1-\lambda)}{\varepsilon}\int_{\Omega}\rho\psi_0|\varphi|^2
\\ \noalign{\medskip}\displaystyle
\phantom{\frac{d}{dt}\int_{\Omega}\rho\psi_0|\varphi|^2}
-2\lambda\int_{\Omega}\rho\psi_1|\varphi|^2-\varepsilon\int_{\Omega}\rho\psi_2|\varphi|^2+(1+\lambda)\int_{\Gamma_0}\rho\psi_0|\varphi|^2
\\ \noalign{\medskip}\displaystyle
 \phantom{\frac{d}{dt}\int_{\Omega}\rho\psi_0|\varphi|^2}
 -\int_{\Gamma_1}\rho((1-\lambda)\psi_0+2\varepsilon\psi_1)|\varphi|^2-\varepsilon\frac{d}{dt}\int_{\Gamma}\rho\psi_0|\varphi|^2.
\end{array}
$$
Observe that, thanks to the choice of the function $\Psi$, we have that $\psi_{0|\Gamma_1}=\psi_{1|\Gamma_1}=0$ and so the sixth term in the right-hand side vanishes. Since $\lambda\in (0,1)$, the second term is positive. Consequently,
$$
\frac{d}{dt}\left(\int_{\Omega}\rho\psi_0|\varphi|^2+\varepsilon\int_{\Gamma}\rho\psi_0|\varphi|^2\right)\geq
-\int_{\Omega}(2\lambda\rho\psi_1+\varepsilon\rho\psi_2)|\varphi|^2.
$$
Since the supports of the functions $\psi_1(t,\cdot)$ and $\psi_2(t,\cdot)$ are included in $\mathbb{R}^{n-1}\times (-\infty,-2\delta)$, we obtain~:
 $$
 \frac{d}{dt}\|(\rho\psi_0)^{1/2}\varphi\|^2_X\geq -Ce^{-2\delta\lambda/\varepsilon}\int_{\Omega}|\varphi|^2.
 $$
 Then, from Remark \ref{IDK}, we deduce that
 $$
 \frac{d}{dt}\|(\rho\psi_0)^{1/2}\varphi\|^2_X\geq-Ce^{-2\delta\lambda/\varepsilon}\|\varphi_T\|^2_X=-Ce^{-2\delta\lambda/\varepsilon}.
 $$
 Integrating between $t$ and $T$, we have :
 $$
 \|(\rho\psi_0(t))^{1/2}\varphi(t)\|^2_X\leq \|(\rho\psi_0(T))^{1/2}\varphi_T\|_X^2+Ce^{-2\delta\lambda/\varepsilon}\leq
 e^{-\delta\lambda/\varepsilon}\|\varphi_T\|_X^2+Ce^{-2\delta\lambda/\varepsilon}\leq Ce^{-\delta\lambda/\varepsilon}.
 $$
 Finally, since $\psi_0(t)_{|\Gamma_0}=\rho_{|\Gamma_0}=1$, we find in particular
 $$
 \varepsilon\|\varphi(t)\|^2_{L^2(\Gamma_0)}\leq Ce^{-\delta\lambda/\varepsilon}\quad t\in (0,T).
 $$
 This gives the desired result (\ref{numberone}).

 \item Proof of (\ref{numbertwo}).

 In this part we prove a quasi-conservation result for the $X$-norm of $\varphi$ (solution of $(S')$ associated to $\varphi_T$) if $\varepsilon$ is small enough.
Let $\theta$ be the solution of the transport equation
$$\left\{ \begin{array}{c}
     \theta_t+\partial_{x_n} \theta = 0\\
\theta(T,.)= \varphi_T \\
      \end{array} \right. \begin{array}{c}
     \text{in } (0,T) \times \Omega,\\
     \text{in } \Omega.
   \end{array}
   $$
One notes that, in fact,
$$ \forall (t,x) \in (0,T) \times \Omega, \ \theta(t,x)=\varphi_T(x',T-t+x_n)$$
and, consequently, thanks to $4\delta < 1-T$,
$$ \theta=\theta_t=\partial_{x_n}\theta=0 \text{ on } (0,T)\times \partial \Omega.$$
We then multiply the equation satisfied by $\varphi$ (see $(S')$) by $\theta $ and we integrate it over $(0,T) \times \Omega$ to get, after integration by parts,
$$
\int_{\Omega} \theta(T,.) \varphi_T - \int_\Omega \theta(0,.) \varphi(0,.)+ \varepsilon \int_{0}^T \int_\Omega \Delta \theta \varphi=0.
$$
Using $\theta(T,.)= \varphi_T$ and Remark \ref{IDK}, one gets for some $C>0$,
$$\|\varphi(0)\|_X \ge \int_\Omega \theta(0,.) \varphi(0,.) \ge 1-C\varepsilon$$
so that, for $\varepsilon< \frac{1}{2C}$,
\begin{equation}\label{equation0}
\|\varphi(0)\|_X \ge \frac{1}{2}.
\end{equation}
This gives (\ref{numbertwo}).

The proof of Theorem \ref{noncontroln} is complete.
\end{itemize}

\appendix
\section{Proof of Proposition \ref{The:Carleman}}
We will use the following notations : $q := (0, T)\times (-1,0)$ , $\sigma:=(0, T)\times \{-1,0\}$, $\sigma_0:=(0, T)\times \{0\}$ and $\sigma_1:=(0, T)\times \{-1\}$. We will now explain how to get the following result.

We perform the proof of this theorem for smooth solutions, so that the general proof follows from a density argument.

We recall the following properties of the weight functions:
\begin{equation}\label{estimationspoids}
\begin{array}{c}
| \alpha_t | \lesssim T \phi^2, \,\,| \alpha_{x t} | \lesssim
     T \phi^2, \,\,| \alpha_{t t} | \lesssim T^2 \phi^3,
\\ \noalign{\medskip}
\alpha_x = -  \phi, \,\,\alpha_{x x} = - \phi
\end{array}
\end{equation}
and we follow the standard method introduced in \cite{FI}. Let $\psi \assign \varphi e^{- s \alpha}$ ; then, using the equation satisfied by $\varphi$, we find
\[
P_1 \psi + P_2 \psi = P_3 \psi \quad\hbox{in }q,
\]
where
\begin{equation}\label{P1}
P_1 \psi = \psi_t + 2 \varepsilon s \alpha_x \psi_x+\psi_x,
\end{equation}
\begin{equation}\label{P2}
P_2 \psi = \varepsilon\psi_{x x} + \varepsilon s^2 \alpha^2_x \psi + s \alpha_t \psi+s\alpha_x\psi - a\psi ,
\end{equation}
and
$$
P_3 \psi= -\varepsilon s \alpha_{x x} \psi.
$$
On the other hand, the boundary conditions are:
\begin{equation}\label{enzero}
  \psi_t + s \alpha_t \psi -\psi_x - s \alpha_x
  \psi - \varepsilon^{-1}\psi = 0\,\,\,\,\text{on }\sigma_0,
\end{equation}
\begin{equation}\label{en-L}
\psi_t + s \alpha_t \psi + \psi_x + s \alpha_x \psi = 0\,\,\,\,\text{on }\sigma_1.
\end{equation}
We take the $L^2$ norm in both sides of the identity in $q$:
\begin{equation}\label{L2}
\|P_1\psi\|^2_{L^2(q)}+\|P_2\psi\|^2_{L^2(q)}+2(P_1\psi,P_2\psi)_{L^2(q)}=\|P_3\psi\|^2_{L^2(q)}.
\end{equation}
Using (\ref{estimationspoids}), we directly obtain
\begin{equation}\label{p3}
\left. \left. \right\| P_3 \psi \right\|^2_{L^2 (q)} \lesssim
    \varepsilon^2 s^2
   \int_q \phi^2 | \psi |^2.
\end{equation}

We focus on the expression of the product $(P_1 \psi, P_2 \psi)_{L^2 (q)}$.
This product contains 15 terms which will be denoted by $T_{i j} (\psi)$ for $1 \leqslant i \leqslant 3$, $1 \leqslant j \leqslant 5$.

\begin{itemize}
\item For the first term in $P_1\psi$, we integrate by parts in time and space. Using that $\psi_{|t=T}=\psi_{|t=0}=0$ and that $a$ is constant, we have
\begin{eqnarray}\label{P11}
\sum_{i=1}^5T_{1i}&=&\int_{q}\psi_t(\varepsilon\psi_{x x} + \varepsilon s^2 \alpha^2_x \psi + s \alpha_t \psi+s\alpha_x\psi - a\psi)\nonumber \\
&=&-\varepsilon s^2\int_{q}\alpha_x\alpha_{xt}|\psi|^2-\frac{s}{2}\int_q (\alpha_{tt}+\alpha_{xt})|\psi|^2+\varepsilon\int_\sigma \psi_t \partial_\nu \psi \nonumber \\
&\gtrsim& -sT(\varepsilon s+T+T^2)\int_q\phi^3|\psi|^2.
\end{eqnarray}
In order to obtain the last estimate, we have used (\ref{estimationspoids}) and the boundary conditions.

\item For the second term in $P_1\psi$, we first have :
\begin{equation}\label{T21}
T_{21}=-\varepsilon^2s\int_{\sigma_0}\phi|\psi_x|^2+\varepsilon^2s\int_{\sigma_1}\phi|\psi_x|^2+\varepsilon^2s\int_{q}\phi|\psi_x|^2.
\end{equation}

Integrating by parts in space, we find
\begin{equation}\label{T22}
T_{22}=-\varepsilon^2s^3\int_{\sigma_0}\phi^3|\psi|^2+\varepsilon^2s^3\int_{\sigma_1}\phi^3|\psi|^2+3\varepsilon^2s^3\int_{q}\phi^3|\psi|^2
\end{equation}
and
\begin{eqnarray}\label{T23}
\sum_{i=3}^5T_{2i}&=&\varepsilon s\int_{\sigma_0}\alpha_x (s\alpha_t+s\alpha_x-a)|\psi|^2-\varepsilon s\int_{\sigma_1}\alpha_x (s\alpha_t+s\alpha_x-a)|\psi|^2\nonumber \\
&-&\varepsilon s\int_{q}[\alpha_{xx} (s\alpha_t+2s\alpha_x-a)+s\alpha_x\alpha_{xt}]|\psi|^2\nonumber \\
&\gtrsim& -\varepsilon Ts^2\int_{\sigma_0}\phi^3|\psi|^2-\varepsilon s[s(T+T^2)+aT^4]\left(\int_{\sigma_1}\phi^3 |\psi|^2
+\int_q\phi^3|\psi|^2\right),
\end{eqnarray}
where we have used estimates (\ref{estimationspoids}).

\item Finally, for the third term in $P_1\psi$ we obtain :
\begin{equation}\label{T31}
\begin{array}{l}\displaystyle
T_{31}+T_{32}=\frac{\varepsilon}{2}\left(\int_{\sigma_0}|\psi_x|^2-\int_{\sigma_1}|\psi_x|^2+s^2\int_{\sigma_0}\phi^2|\psi|^2
-s^2\int_{\sigma_1}\phi^2|\psi|^2\right)-\varepsilon s^2\int_q\phi^2|\psi|^2
\\ \noalign{\medskip}\displaystyle
\phantom{T_{31}+T_{32}}\gtrsim -\varepsilon T^2\left(\int_{\sigma_1}\phi|\psi_x|^2+s^2\int_{\sigma_1}\phi^3|\psi|^2+s^2\int_{q}\phi^3|\psi|^2\right)
\end{array}
\end{equation}
and
\begin{equation}\label{T33}
\begin{array}{l}\displaystyle
\sum_{i=3}^5T_{3i}=\frac{1}{2}\left(\int_{\sigma_0}(s\alpha_t+s\alpha_x-a)|\psi|^2-\int_{\sigma_1}(s\alpha_t+s\alpha_x-a)|\psi|^2
-s\int_q(\alpha_{tx}+\alpha_{xx})|\psi|^2\right)
\\ \noalign{\medskip}\displaystyle
\phantom{\sum_{i=3}^5T_{3i}}\gtrsim -s(T^3+T^4)\left(\int_{\sigma_0}\phi^3|\psi|^2+\int_{\sigma_1}\phi^3|\psi|^2+\int_{q}\phi^3|\psi|^2\right)
-aT^6\int_{\sigma_0}\phi^3|\psi|^2.
\end{array}
\end{equation}
\end{itemize}
Putting together (\ref{P11})-(\ref{T33}), we obtain, since $T^2 \lesssim T+T^3$,
\begin{eqnarray}\label{produitdouble}
(P_1\psi,P_2\psi)_{L^2(q)}=\sum\limits_{\substack{ 1 \le i \le 3 \\ 1 \le j \le 5}} T_{ij} (\psi) &\ge& \varepsilon^2s^3\left(\int_{\sigma_1}\phi^3|\psi|^2+\int_q\phi^3|\psi|^2\right)
+\varepsilon^2s\left(\int_{\sigma_1}\phi|\psi_x|^2+\int_q\phi|\psi_x|^2\right)\nonumber\\
&-& C\left(sT[\varepsilon(sT+s+aT^3)+T+T^3]\left(\int_{\sigma_1}\phi^3|\psi|^2+\int_q\phi^3|\psi|^2\right) \right. \nonumber \\
&+& \left. [s^2\varepsilon(T+\varepsilon s)+T^3(s+sT+aT^3)]\int_{\sigma_0}\phi^3|\psi|^2 \right. \nonumber\\
&+&\left. \varepsilon^2s\int_{\sigma_0}\phi|\psi_x|^2 +\varepsilon T^2\int_{\sigma_1}\phi|\psi_x|^2 \right).
\end{eqnarray}
We readily observe that the second line of this expression can be absorbed by the first term in the right-hand side of the first line, that is to say,
$$
\varepsilon^2s^3\left(\int_{\sigma_1}\phi^3|\psi|^2+\int_q\phi^3|\psi|^2\right),
$$
provided that
\begin{equation}\label{choixepsilon}
s \gtrsim \varepsilon^{-1} (T + T^2)+a^{1/2} \varepsilon^{-1/2} T^2.
\end{equation}
Consequently, we obtain
\begin{eqnarray}\label{produitdouble2}
(P_1\psi,P_2\psi)_{L^2(q)} &\ge& \frac{\varepsilon^2s^3}{2}\left(\int_{\sigma_1}\phi^3|\psi|^2+\int_q\phi^3|\psi|^2\right)
+\varepsilon^2s\left(\int_{\sigma_1}\phi|\psi_x|^2+\int_q\phi|\psi_x|^2\right)\nonumber \\
&-&C\left([s^2\varepsilon(T+\varepsilon s)+T^3(s+sT+aT^3)]\int_{\sigma_0}\phi^3|\psi|^2+\varepsilon^2s\int_{\sigma_0}\phi|\psi_x|^2 \right. \nonumber\\
&+& \left. \varepsilon T^2\int_{\sigma_1}\phi|\psi_x|^2 \right).
\end{eqnarray}
Furthermore, the last term in this expression is absorbed by
$$
\varepsilon^2s\int_{\sigma_1}\phi^3|\psi|^2
$$
if $s\gtrsim  \varepsilon^{-1}T^2$. We also observe that the term in $\sigma_0$ can be estimated as follows:
$$
[s^2\varepsilon(T+\varepsilon s)+T^3(s+sT+aT^3)]\int_{\sigma_0}\phi^3|\psi|^2\lesssim \varepsilon^2s^3\int_{\sigma_0}\phi^3|\psi|^2,
$$
provided that $s \gtrsim T^2(\varepsilon^{-1}+a^{1/3}\varepsilon^{-2/3})$. This choice of the parameter $s$ is implied by (\ref{choixepsilon}).

Coming back to (\ref{L2}), we have proved that
\begin{equation}\label{produitdouble2}
\begin{array}{l}\displaystyle
\|P_1\psi\|^2_{L^2(q)}+\|P_2\psi\|^2_{L^2(q)}+ \varepsilon^2s^3\left(\int_{\sigma_1}\phi^3|\psi|^2+\int_q\phi^3|\psi|^2\right)
+\varepsilon^2s\left(\int_{\sigma_1}\phi|\psi_x|^2+\int_q\phi|\psi_x|^2\right)
\\ \noalign{\medskip}\displaystyle
\lesssim \varepsilon^2s^3\int_{\sigma_0}\phi^3|\psi|^2+\|P_3\psi\|^2_{L^2(q)}+\varepsilon^2s\int_{\sigma_0}\phi|\psi_x|^2.
\end{array}
\end{equation}
for $s$   as in (\ref{choixepsilon}). Observe that from (\ref{P1}) and $s \varepsilon \gtrsim T^2$, we deduce that
\begin{equation}\label{produitdouble22}
\begin{array}{l}\displaystyle
s^{-1}\int_q\phi^{-1}|\psi_t|^2+\varepsilon^2s^3\left(\int_{\sigma_1}\phi^3|\psi|^2+\int_q\phi^3|\psi|^2\right)
+\varepsilon^2s\left(\int_{\sigma_1}\phi|\psi_x|^2+\int_q\phi|\psi_x|^2\right)
\\ \noalign{\medskip}\displaystyle
\lesssim \varepsilon^2s^3\int_{\sigma_0}\phi^3|\psi|^2+\|P_3\psi\|^2_{L^2(q)}+\varepsilon^2s\int_{\sigma_0}\phi|\psi_x|^2.
\end{array}
\end{equation}
The term in $P_3\psi$ can be absorbed by the term in the left-hand side thanks to (\ref{p3}) and for $s \gtrsim \varepsilon^{-1}T^2 \gtrsim  T^2$.
We finally estimate
$$
\varepsilon^2s\int_{\sigma_0}\phi|\psi_x|^2
$$
using the boundary condition at $x=0$ given by (\ref{enzero}). It follows that
$$
\varepsilon^2s\int_{\sigma_0}\phi|\psi_x|^2\lesssim\varepsilon^2s\left(\int_{\sigma_0}\phi(s^2(\alpha_t)^2+s^2(\alpha_x)^2+\varepsilon^{-2})|\psi|^2
+\int_{\sigma_0}\phi|\psi_t|^2\right)
$$
Using (\ref{estimationspoids}), we find
\begin{equation}\label{Betis2}
\varepsilon^2s\int_{\sigma_0}\phi|\psi_x|^2 \lesssim \varepsilon^2s^3\int_{\sigma_0}(\phi^3+T^2\phi^5)|\psi|^2
+\varepsilon^2s\int_{\sigma_0}\phi|\psi_t|^2.
\end{equation}

\noindent We now come back to $\varphi$, recalling that $\psi=e^{-s\alpha}\varphi$. Then, using again (\ref{estimationspoids}) and \eqref{choixepsilon}, we get from (\ref{produitdouble22}) and (\ref{Betis2})
\begin{equation}\label{Seville}
\begin{array}{l}\displaystyle
s^{-1}\int_q\phi^{-1}e^{-2s\alpha}|\varphi_t|^2+\varepsilon^2s^3\left(\int_q\phi^3e^{-2s\alpha}|\varphi|^2+\int_{\sigma}\phi^3e^{-2s\alpha}|\varphi|^2\right)
+\varepsilon^2s\int_{\sigma_1}\phi e^{-2s\alpha}|\varphi_x|^2
\\ \noalign{\medskip}\displaystyle
\lesssim \varepsilon^2s^3\int_{\sigma_0}(\phi^3+T^2\phi^5)e^{-2s\alpha}|\varphi|^2+
\varepsilon^2s\int_{\sigma_0}\phi e^{-2s\alpha}|\varphi_t|^2.
\end{array}
\end{equation}
The last step is to estimate the term in $|\varphi_t|^2$ on $\sigma_0$ in the right-hand side of (\ref{Seville}). Using that $s\varphi \gtrsim 1$ for $s$ satisfying \eqref{choixepsilon}, we have
\begin{eqnarray}\label{Heliopolis}
\varepsilon^2s\int_{\sigma_0}\phi e^{-2s\alpha}|\varphi_t|^2&=&-\varepsilon^2s\int_{\sigma_0}\phi e^{-2s\alpha}\varphi_{tt}\varphi+\frac{\varepsilon^2s}{2}\int_{\sigma_0}(\phi e^{-2s\alpha})_{tt}|\varphi|^2\nonumber\\
&\lesssim& \varepsilon^2s\int_{\sigma_0}\phi e^{-2s\alpha}|\varphi_{tt}||\varphi|+\varepsilon^2T^2s^3\int_{\sigma_0}\phi^5 e^{-2s\alpha}|\varphi|^2.
\end{eqnarray}
\par The goal is now to estimate $\varphi_{tt}$ on $\sigma_0$. For this purpose, let us set $\rho(t):= s^{-5/2}\phi(t,-1)^{-5/2}e^{-s\alpha(t,-1)}$ and
$w^*:=\rho\varphi_t$. Then, $w^*$ satisfies
$$
 (S'^a_*) \left\{ \begin{array}{c}
     w^*_t + w^*_x + \varepsilon w^*_{xx} -a w^* = \rho'\varphi_t\\
     \varepsilon (w^*_t - \partial_{\nu} w^*) - w^* = \varepsilon\rho'\varphi_t\\
     w^*_t - \partial_{\nu} w^* = \rho'\varphi_t\\
     w^* (T, .) = 0
   \end{array} \right. \begin{array}{c}
     \text{in } (0, T) \times (-1,0),\\
     \text{on } (0, T) \times \{0\},\\
     \text{on } (0, T) \times \{-1\},\\
     \text{in } (-1,0).
   \end{array}
$$
$\bullet$ In a first step, we multiply this system by $w^*$ and we integrate in $q$. After some computations, we obtain
$$
\varepsilon\int_q(w^*_x)^2+\frac{1}{2}\int_{\sigma}|w^*|^2+a\int_q|w^*|^2=\int_q\rho'\varphi_tw^*-\varepsilon\int_{\sigma}\rho'\varphi_tw^*.
$$
In particular, we have, using Young inequality,
\begin{equation}\label{Lyon}
\varepsilon\int_q(w^*_x)^2 \lesssim \int_q|\rho'\varphi_t|^2+\varepsilon^2\int_{\sigma}|\rho'\varphi_t|^2.
\end{equation}
$\bullet$ Then, we multiply by $\varepsilon w^*_t$. Analogously, we get
$$
\frac{\varepsilon}{2}\int_q(w^*_t)^2+\frac{\varepsilon}{2}\int_{\sigma}(w^*_t)^2 \lesssim \varepsilon^2\int_{\sigma}|\rho'\varphi_t|^2+
\varepsilon\int_{q}|\rho'\varphi_t|^2+\varepsilon\int_q(w^*_x)^2.
$$
Combining this with (\ref{Lyon}), we obtain
$$
\varepsilon^2\int_{\sigma}|w_t|^2 \lesssim \varepsilon^2\int_{\sigma}|\rho'\varphi_t|^2+\int_q|\rho'\varphi_t|^2.
$$
Since $w_t=\rho'\varphi_t+\rho\varphi_{tt}$, we have
$$
\varepsilon^2\int_{\sigma}\rho^2|\varphi_{tt}|^2 \lesssim \varepsilon^2\int_{\sigma}|\rho'\varphi_t|^2+\int_q|\rho'\varphi_t|^2.
$$
In particular, we find
\begin{equation}\label{Marseille}
\varepsilon^2s^{-5}\int_{\sigma_0}\phi^{-5}(t,-1)e^{-2s\alpha(t,-1)}|\varphi_{tt}|^2 \lesssim \varepsilon^2 s^{-1}\int_{\sigma}\phi^{-1}e^{-2s\alpha}|\varphi_t|^2+s^{-1}\int_q\phi^{-1}e^{-2s\alpha}|\varphi_t|^2.
\end{equation}
Here, we have used that
$$
\phi^{-1}(t,-1)e^{-2s\alpha(t,-1)}\leq \phi^{-1}(t,x)e^{-2s\alpha(t,x)}\quad\hbox{ for all }x\in (-1,0).
$$
\vskip0.2cm \noindent Coming back to (\ref{Heliopolis}), we have
$$
\begin{array}{l}\displaystyle
\varepsilon^2s\int_{\sigma_0}\phi e^{-2s\alpha}|\varphi_t|^2\leq C\varepsilon^2 s^7\int_{\sigma_0}\phi^7e^{-4s\alpha+2s\alpha(t,-1)}|\varphi|^2
\\ \noalign{\medskip}\displaystyle
\phantom{\varepsilon^2s\int_{\sigma_0}\phi(t,-1) e^{-2s\alpha(t,-1)}(\varphi_t)^2}+\delta \varepsilon^2s^{-5}\int_{\sigma_0}\phi^{-5}(t,-1)e^{-2s\alpha(t,-1)}|\varphi_{tt}|^2,
\end{array}
$$
for $s \gtrsim T^{2}$ and all $\delta>0$. From (\ref{Marseille}), we now obtain
$$
\begin{array}{l}\displaystyle
\varepsilon^2s\int_{\sigma_0}\phi e^{-2s\alpha}|\varphi_t|^2\leq C\varepsilon^2s^7\int_{\sigma_0}\phi^7e^{-4s\alpha+2s\alpha(t,-1)}|\varphi|^2
\\ \noalign{\medskip}\displaystyle
+C \delta\left(\varepsilon^2 s^{-1}\int_{\sigma}\phi^{-1}e^{-2s\alpha}|\varphi_t|^2+s^{-1}\int_q\phi^{-1}e^{-2s\alpha}|\varphi_t|^2\right).
\end{array}
$$
Combining this with (\ref{Seville}), using the boundary conditions and taking $\delta$ small enough, we conclude that if $s$ satisfies \eqref{choixepsilon},
$$
\begin{array}{l}\displaystyle
s^{-1}\int_q\phi^{-1}e^{-2s\alpha}|\varphi_t|^2+\varepsilon^2s^3\left(\int_q\phi^3e^{-2s\alpha}|\varphi|^2+\int_{\sigma}\phi^3e^{-2s\alpha}|\varphi|^2\right)
+\varepsilon^2s\int_{\sigma_1}\phi e^{-2s\alpha}|\varphi_x|^2
\\ \noalign{\medskip}\displaystyle
\lesssim
\varepsilon^2s^7\int_{\sigma_0}\phi^7 e^{-4s\alpha+2s\alpha(t,-1)}|\varphi|^2.
\end{array}
$$
In particular, this implies the desired inequality (\ref{Carlemana}).

\end{document}